\documentclass[12pt,a4paper]{amsart}
\usepackage{amsfonts}
\usepackage{ifthen}
\usepackage{amsmath}
\usepackage{graphicx}
\usepackage{amscd,amssymb,amsthm}

\newcounter{minutes}
\divide\time by 60
\newcounter{hours}
\multiply\time by 60 \addtocounter{minutes}{-\time}

\title{\bf On quotients and differences of hypergeometric functions}

\author{ Slavko Simi\'c}
\author{Matti Vuorinen}

\address{ Mathematical Institute SANU, Kneza Mihaila 36, 11000
Belgrade, Serbia} \email{ ssimic@turing.mi.sanu.ac.rs}
\address{Department of Mathematics, University of Turku, 20014 Turku,
Finland} \email{vuorinen@utu.fi}

\newtheorem{theorem}[equation]{Theorem}

\newtheorem{lemma}[equation]{Lemma}

\newtheorem{corollary}[equation]{Corollary}
\newtheorem{remark}[equation]{Remark}

\newtheorem{question}[equation]{Question}
\newtheorem{nonsec}[equation]{}

\numberwithin{equation}{section}

\begin{document}

\begin{abstract}
For Gaussian hypergeometric functions $F(x)= F(a,b;c;x),$ $a,b,c>0,$
we consider the quotient
$ Q_F(x,y)= (F(x)+F(y))/F(z)$ and the difference $ D_F(x,y)= F(x)+F(y)-F(z)$
for $0<x,y<1$ with $z=x+y-xy \,.$ We give best possible bounds for  both expressions
under various hypotheses about the parameter triple $(a,b;c)\,.$
\end{abstract}

\def\thefootnote{}
\footnotetext{ \texttt{\tiny File:~\jobname .tex,
          printed: \number\year-\number\month-\number\day,
          \thehours.\ifnum\theminutes<10{0}\fi\theminutes}
} \makeatletter\def\thefootnote{\@arabic\c@footnote}\makeatother

\maketitle

{\small \sc Keywords.}{ Sub-additivity; Hypergeometric functions;
Inequalities.}

{\small \sc 2010 Mathematics Subject Classification.}{   26D06,33C05}

\section{Introduction}

Among special functions, the hypergeometric function has perhaps
the widest range of applications. For instance, several well-known classes
of special functions such as complete elliptic integrals, Legendre functions,
Chebyshev and Jacobi polynomials, and some elementary functions, such as
the logarithm, are particular cases of it, cf. \cite{as}.
In a recent paper \cite{kmsv} the authors studied various extensions of the
Bernoulli inequality for functions of logarithmic type. In
particular, the  zero-balanced hypergeometric function
$F(a,b;a+b;x), a,b>0,$ occurs in these studies, because it has a
logarithmic singularity at $x=1\,,$ see \eqref{hypergeom} below.
We now continue the discussion of some of the questions for
quotients and differences of hypergeometric functions that were
left open in \cite{kmsv}.

Motivated by the asymptotic behavior of the function $F(x)=F(a,b;c;x)$
when $x\to 1^-$, see \eqref{hypergeom} and \eqref{gIdty}, we define for $0<x,y<1$
$$
Q_F(x,y):=\frac{F(x)+F(y)}{F(x+y-xy)} \,, \quad
D_F(x,y):=F(x)+F(y)-F(x+y-xy) \,.
$$

Our task in this paper is to give tight bounds for these quotients and
differences assuming various relationships between the parameters $a,b,c$.

For the general case we can formulate the next

\begin{theorem} \label{quotF}
For $a,b,c>0$ and $0<x,y<1$ let
\begin{equation} \label{quotF1}
Q_F(x,y)= \dfrac{F(a,b;c;x)+F(a,b;c;y)}{F(a,b;c;x+y-xy)}\,.
\end{equation}
Then
\begin{equation} \label{quotF2}
0 < Q_F(x,y)\le 2.
\end{equation}
\end{theorem}

The bounds in \eqref{quotF2} are best possible as can be seen by taking
\cite[15.1.8]{as}
$$
F(x)=F(a,b;b;x)=(1-x)^{-a}:=F_0(x).
$$

Then
$$
Q_{F_0}(x,y)=\frac{(1-x)^{-a}+(1-y)^{-a}}{(1-x-y+xy)^{-a}}=
\frac{(1-x)^{-a}+(1-y)^{-a}}{((1-x)(1-y))^{-a}}=(1-x)^a+(1-y)^a,
$$
and the conclusion follows immediately.
Similarly,

\begin{theorem} \label{difF}
For $a,b>0, c>a+b$ and $0<x,y<1$, we have
\begin{equation} \label{difF1}
|D_F(x,y)|\le A,
\end{equation}
with $A=A(a,b,c)=\frac{\Gamma(c)\Gamma(c-a-b)}{\Gamma(c-a)\Gamma(c-b)}=F(a,b;c;1)$
as the best possible constant.
\end{theorem}

\bigskip

Most intriguing is the zero-balanced case. For example,

\begin{theorem} \label{diffF}
For $a,b>0$ and $0<x,y<1$ let
\begin{equation} \label{diffF1}
D_F(x,y)= F(a,b;a+b;x)+F(a,b;a+b;y) -F(a,b;a+b;x+y-xy) \,.
\end{equation}
Then
\begin{equation} \label{diffF2}
\dfrac{2R}{B} -1 < D_F(x,y)\le 1 \,
\end{equation}
with $R=R(a,b)= -2\gamma -\psi(a)-\psi(b), B=B(a,b)\,.$
\end{theorem}

The constant $1$ on the right-hand side is the best possible upper bound,
but  the best constant on the left-hand side is not known to us.

The proofs of these results are based on the monotonicity of some
functions and thus they hold under much weaker hypotheses than
stated above, see Remark \ref{moregeneral}.

\medskip

In the sequel we shall give a complete answer to an open question
posed in \cite{kmsv}.


\section{Preliminary results}

In this  section we recall some well-known properties of the
Gaussian hypergeometric function $F(a,b;c;x)$ and certain of its
combinations with other functions, for further applications.
For basic information, the handbooks \cite{as,olbc} may be recommended.

It is well known that hypergeometric functions are closely
related to the classical {\it gamma function} $\Gamma(x)$,
the {\it psi function} $\psi(x)$, and the {\it beta function}
$B(x,y)$.
For Re$\thinspace x>0,~~$ Re$\thinspace y>0$, these functions
are defined by

\begin{equation} \label{myB}
\Gamma(x)\equiv \int\limits_0^{\infty}e^{-t}t^{x-1}dt, ~
\psi(x)\equiv \frac{\Gamma'(x)}{ \Gamma(x)}, ~
B(x,y)\equiv {\frac{\Gamma(x)\Gamma(y)}{\Gamma(x+y)}},
\end{equation}
respectively (cf. \cite[Chap. 6]{as}). It is well known that the gamma
function
satisfies the {\it difference equation} \cite[6.1.15]{as}

\begin{equation}
\Gamma(x+1)=x\Gamma(x),
\end{equation}
and the {\it reflection property} \cite[6.1.17]{as}

\begin{equation}
\Gamma(x)\Gamma(1-x)=\displaystyle\frac{\pi}{  {\sin \pi x}}=B(x,1-x).
\end{equation}
We shall also need the function
\begin{equation} \label{myR}
R(a,b)\! \equiv\! -2\gamma\!-\!\psi(a)\!-\!\psi(b),~ R(a)\! \equiv\!
R(a,1\!-\! a), ~ R\bigl(\frac{1}{2}\bigr)\! =\! \log 16,
\end{equation}
where $\gamma$ is the {\it Euler-Mascheroni constant} given
by
\begin{equation}
\gamma=\lim\limits_{n\to \infty}\biggl (\sum_{k=1}^{n}
{1\over k}-\log n \biggr ) = 0.577215\ldots.
\end{equation}


Given complex numbers
$a,b,$ and $c$ with $c\neq 0,-1,-2, \dots $,
the {\em Gaussian hypergeometric function} is the analytic
continuation to the slit plane ${\mathbb C} \setminus [1,\infty)$ of
the series
\begin{equation} \label{eq:hypdef}
F(a,b;c;z) = {}_2 F_1(a,b;c;z) =
\sum_{n=0}^{\infty} \frac{(a,n)(b,n)}{(c,n)} \frac{z^n}{n!}\,,\:\: |z|<1 \,.
\end{equation}
Here $(a,0)=1$ for $a \neq 0$, and $(a,n)$
is the {\em shifted factorial function} or the {\em Appell symbol}
$$
(a,n) = a(a+1)(a+2) \cdots (a+n-1)
$$
for $n \in {\mathbb N} \setminus \{0\}$, where
$ {\mathbb N} = \{ 0,1,2,\ldots\}$.

The hypergeometric function has the following
simple differentiation formula (\cite[15.2.1]{as})
\begin{equation}
{d\over dx}F(a,b;c;x)={ab\over c}F(a+1,b+1;c+1;x) \,.
\end{equation}

An  important tool for our work is the following classification of the
behavior of the hypergeometric function near $x=1$ in the three cases
$ a+b < c,
a+b = c,$ and $ a+b > c :$
\begin{equation}\label{hypergeom}
\left\{ \begin{array}{l}
F(a,b;c;1) = \dfrac{\Gamma(c) \Gamma(c-a-b)}{\Gamma(c-a) \Gamma(c-b)}, ~~ a+b <
c,\\[.5cm]
B(a,b)F(a,b;a\!+\!b;x)\!+\!\log(1\!-\!x) \!=\! R(a,b)\!+\!{\rm O}
 ((1\!-\!x)\log(1\!-\!x)), \\[.4cm]
 F(a,b;c;x) = (1\!-\!x)^{c -a -b} F(c\!-\!a,c\!-\!b;c;x), ~ c < a\!+\!b.
\end{array} \right.
\end{equation}

Some basic properties of this series may be found in standard
handbooks, see for example \cite{as}. For some rational triples $(a,b,c)$,
the function $F(a,b;c;x)$ can be expressed in terms of well-known elementary
function. A particular case that is often used in this paper is \cite[15.1.3]{as}
\begin{equation}\label{functiong}
  g(x) \equiv x F(1,1;2;x) = \log \frac{1}{1-x}.
\end{equation}
It is clear that for $a,b,c>0$ the function $F(a,b;c;x)$ is a
strictly increasing map from  $[0,1)$ into $[1,\infty)\,.$ For
$a,b>0$ we see by (\ref{hypergeom}) that $xF(a,b;a+b;x)$ defines
an increasing homeomorphism from $[0,1)$ onto $[0,\infty)\,.$

\begin {theorem}{\cite{abrvv},\cite[Theorem 1.52]{avv}} \label{1.57}
For $a, b > 0,$ let $B=B(a,b)$ be as in $\eqref {myB}$, and let $R=R(a,b)$
be as in $\eqref{myR}$. Then the following are true.

\begin {itemize}
\item[(1)] The function $f_1(x) \equiv \frac{ F(a,b;a+b;x)-1}{
\log (1/ (1-x))} $ is strictly increasing from $(0,1)$ onto
\smallbreak
\noindent
$(ab/(a+b),1/B).$

\item[(2)] The function $f_2(x) \equiv BF(a,b;a+b;x) + \log (1-x)$
is strictly decreasing from $(0,1)$ onto $(R,B).$

\item[(3)] The function $f_3(x) \equiv BF(a,b;a+b;x) + (1/x)\log (1-
x)$ is increasing from $(0,1)$ onto $(B-1,R)$ if $a,b \in
(0,1).$

\item[(4)] The function $f_3$ is decreasing from $(0,1)$ onto
$(R,B-1)$ if $a,b \in (1,\infty ).$

\item[(5)] The function
$$f_4(x) \equiv \frac{xF(a,b;a+b;x)}  {\log (1/ (1-x))}$$

is decreasing from $(0,1)$ onto $(1/B,1)$ if $a,b \in (0,1)$.

\item[(6)] If $a, b > 1,$ then $f_4$ is increasing from $(0,1)$
onto $(1,1/B).$

\item[(7)] If $a = b = 1$, then $f_4(x) = 1$ for all $x \in
(0,1).$
\end {itemize}
\end {theorem}

\bigskip
We also need the following
refinement of some parts of Theorem \ref{1.57}.

\begin{lemma} \cite[Cor. 2.14]{pv} \label{pvlem}
 For $c,d>0$,  let $B=B(c,d)$ be as in $\eqref {myB}$, and let $R=R(c,d)$
 be as in $\eqref{myR}$ and denote
$$f(x) \equiv \frac{xF(c,d;c+d;x)}  {\log (1/ (1-x))}$$.

(1) If $c\in (0,\infty)$ and $d\in (0,1/c]$, then the function $f$
is decreasing with range $(1/B,1)$;

(2) If $c\in (1/2,\infty)$ and $d\ge c/(2c-1)$, then $f$ is
increasing from $(0,1)$ to the range $(1,1/B)$.

(3) If $c\in (0,\infty)$ and $d\in (0,1/c]$, then the function $h$
defined by
$$
h(x):=BF(c,d;c+d;x)+(1/x)\log(1-x)
$$
is increasing from $(0,1)$ onto $(B-1,R)$.

(4) If $c\in (1/3,\infty)$ and $d\ge (1+c)/(3c-1)$, then $h$ is
increasing from $(0,1)$ onto $(R,B-1)$.
\end{lemma}

For brevity we write $\mathbb R_+ = (0,\infty)\,.$

\begin{lemma} (Cf. \cite[1.24, 7.42(1)]{avv}) \label{mypropo}
(1) If $E(t)/t$ is an increasing function on $\mathbb R_+$, then
$E$ is sub-additive, i.e. for each $x,y>0$ we have that
$$
E(x)+E(y)\le E(x+y)\,.
$$

(2)  If $E(t)/t$ decreases on $\mathbb R_+$, then $E$ is a
super-additive function, that is
$$
E(x)+E(y)\ge E(x+y)
$$
for $x,y\in \mathbb R_+$.
\end{lemma}

\section{Main results}

\bigskip


By \eqref{hypergeom} the zero-balanced hypergeometric function
$F(a,b;a+b;x)$ has a logarithmic singularity at $x=1\,.$
We shall now demonstrate that its behavior is nearly logarithmic
also in the sense that some basic
identities of the logarithm yield functional inequalities for the
zero-balanced function.

Next, writing the basic addition formula for the logarithm
\[
  \log z+\log w =\log(z w), \quad z,w >0 \,,
\]
in terms of the function $g$ in \eqref{functiong}, we have
\begin{equation} \label{gIdty}
  g(x)+g(y) = g(x+y-xy), \quad x,y \in (0,1).
\end{equation}
This identity \eqref{gIdty} shows that $Q_h(x,y)\equiv 1$ and 
$D_h(x,y) \equiv 0$ for $(a,b,c)=(1,1,2)$ and $h(x) = xF(a,b;c;x)\,.$
Based on this observation and a few computer experiments
we posed in \cite{kmsv} the following question:

\begin{question} \label{myq3}
{\rm
 Fix $c,d >0$ and let  $g(x) = x F(c,d;c+d;x)$  for $x \in
(0,1)$ and set
\[
  Q_g(x,y) = \frac{g(x)+g(y)}{g(x+y-xy)}
\]
for $x,y \in (0,1)$.

(1) \ For which values of $c$ and $d\,,$ this function is bounded
from below and above?

(2) Is it true that
\begin{enumerate}
  \item[a)] $Q_g(x,y) \ge 1$, if $cd \le 1$?
  \item[b)] $Q_g(x,y) \le 1$, if $c,d > 1$?
\item[c)] Are there counterparts of Theorem \ref{diffF} for the function
$$  D_g(x,y)= g(x)+g(y)-g(x+y-xy)\,?$$
\end{enumerate}
}
\end{question}

\bigskip

We shall give a complete answer to this question in the sequel.

Note firstly that the quotient $Q_g$ is always bounded. Namely,

\begin{theorem} \label{hRthm}
For all $c,d>0$ and all $x,y\in (0,1)$ we have that
\[
0<Q_g(x,y)<2.
\]
\end{theorem}

 A refinement of these bounds for some particular $(c,d)$ pairs is given by the following two
assertions.

\begin{theorem} \label{hBthm}
(1) \   For $c,d>0, cd\le 1$ and $x,y \in (0,1)$ we have
  \[
    \frac{1}{B(c,d)} \le Q_g(x,y) \le B(c,d).
  \]

(2) \ For $c,d>0, 1/c+1/d \le 2$ and $x,y \in (0,1)$ we have
  \[
    B(c,d) \le Q_g(x,y) \le \frac{1}{B(c,d)}.
  \]
\end{theorem}

Note that parts (1) and (3) of Lemma \ref{pvlem} imply that for $c,d>0,
cd\le 1, (c,d)\neq (1,1)$ we have
\begin{equation} \label{Bgt1}
R(c,d)>0, B(c,d)> 1\, .
\end{equation}

\begin{theorem}
  For $cd \le 1$ and $x,y \in (0,1)$ we have
  \[
    \frac{ B(c,d)-1 }{ R(c,d) } \le Q_g(x,y) \le \frac{ 2 R(c,d) }{ B(c,d)-1 }.
  \]

\end{theorem}

\bigskip

We shall prove now the hypothesis from the second part of Question
\ref{myq3} under the condition $1/c+1/d\le 2\,$ in part $b)$
which, in particular, includes the case $c>1, d>1\,.$

\begin{theorem} \label{ssthm3}
Fix $c,d>0$ and let $Q$ and $g$ be as in Question \ref{myq3}.
\begin{enumerate}
  \item If $cd \le 1$ then $Q_g(x,y) \ge 1\,$ for all $x,y\in
  (0,1)$.

  \item If $1/c+1/d\le 2$, then $Q_g(x,y) \le 1\,$ for all $x,y\in
  (0,1)$.
\end{enumerate}
\end{theorem}

\bigskip

Counterparts of Theorem \ref{diffF} for the difference $D_g$ are given in the
next assertion.

\begin{theorem} \label{ssthm4}
Fix $c,d>0$ and let $D$ be as in Question \ref{myq3}.
\begin{enumerate}
  \item If $cd \le 1$, then
  $$
  0\le D_g(x,y)< \frac{2R(c,d)+1}{B(c,d)}-1 \,
  $$
  for all $x,y\in(0,1)$.

  \item If $1/c+1/d\le 2$, then
  $$
  \frac{2R(c,d)+1}{B(c,d)}-1<D_g(x,y) \le 0\,
  $$

  for all $x,y\in (0,1)$.
\end{enumerate}
\end{theorem}

\bigskip

Combining results above, we obtain the following two-sided bounds
for the quotient $Q_g$.

\bigskip


 \begin{corollary} \label{cor38} Fix $c,d>0$ and let $Q$ be as in Question
\ref{myq3}.
\begin{enumerate}
\item If $cd \le 1$, then
$$ 1\le Q_g(x,y) < \min\{B(c,d), 2\} $$
for all $x,y\in(0,1)$.
\bigskip
\item If $1/c+1/d\le 2$, then
 $$ B(c,d)<Q_g(x,y) \le 1 $$
for all $x,y\in(0,1)$.
\end{enumerate}
\end{corollary}

\bigskip


The assertions above represent a valuable tool for estimating
quotients and differences of a hypergeometric function with different
arguments. To illustrate this point, we give an example.

In the paper \cite{kmsv}, motivated by the relation $g(x)=2g(1-\sqrt{1-x})$ with $g$ as in
\eqref{functiong}, the authors asked the question about the bounds for the function $S(t)$
defined by
$$
S(t):=\frac{g(t)}{g(1-\sqrt{1-t})}, \  t\in(0,1),
$$
where $g(t):=tF(a,b;a+b;t), \ a,b>0$.

\bigskip

An answer follows instantly applying Corollary \ref{cor38}.

 \begin{theorem} \label{ssthm5} \ Let $S(t):=\frac{g(t)}{g(1-\sqrt{1-t})}, \  t\in(0,1),$
 with $g(t):=tF(a,b;a+b;t), \ a,b>0$.

\begin{enumerate}
  \item If $ab \le 1$, then
  $$
  1< S(t) \le 2.
  $$

\bigskip

\item If $1/a+1/b\le 2$, then
 $$
 2\le S(t) < \frac{2}{B(a,b)}.
 $$
\end{enumerate}
\end{theorem}

\bigskip

\section{Proofs}

\bigskip

\begin{nonsec}
{\bf Proof of Theorem \ref{quotF}.}  {\rm
The proof is based solely on the monotonicity property of the function
$F(x)= F(a,b;c;x)$. Namely, for $x,y\in (0,1)$, put $z=x+y-xy, z\in (0,1)$. Since
$$
x\le\max\{x,y\}, y\le\max\{x,y\}; z\ge\max\{x,y\},
$$
and $F(u)$ is monotone increasing in $u$, we obtain
$$
Q_F(x,y)=\frac{F(x)+F(y)}{F(z)}\le
\frac{2F(\max\{x,y\})}{F(\max\{x,y\})}=2.
$$

The left-hand bound is trivial.} \hfill $\square$
\end{nonsec}

\bigskip

\begin{nonsec}
{\bf Proof of  Theorem \ref{difF}.} {\rm The assertion of this
theorem is a consequence of the previous one and
\eqref{hypergeom}. Indeed, from (1.3) we get
$$
-F(a,b;c;z)<D_F(x,y)\le F(a,b;c;z),
$$
that is,}
$$
|D_F(x.y)|\le F(a,b;c;z)=F(a,b;c;1-(1-x)(1-y))\le F(a,b;c;1)=A.
$$
\hfill $\square$
\end{nonsec}

 \begin{nonsec}{\bf Proof of Theorem \ref{diffF}.} \ {\rm By Theorem \ref{1.57}(1) we have that
$$
f_1(x)= \dfrac{F(a,b;a+b;x)-1}{\log(1/(1-x))} =
\dfrac{F(x)-1}{\log(1/(1-x))}
$$
is strictly increasing on $(0,1)\,.$ Putting $x=1-e^{-t}, t>0,$ we
see that
$$
\dfrac{F(1-e^{-t})-1}{t}
$$
is strictly increasing on $(0,\infty)\,,$ that is, by Lemma
\ref{mypropo} the function $G(t)=F(1-e^{-t})-1$ is sub-additive,
i.e.
$$G(u)+G(v) \le G(u+v)\,.$$
Therefore
$$ F(1-e^{-u})+ F(1-e^{-v})- F(1-e^{-(u+v)})\le 1\,. $$
Putting $1-e^{-u}=x,1-e^{-v}=y ,$ we see that
$$ D_F(x,y)=F(x)+F(y)-F(x+y-xy)\le 1\,,$$
and the first part of Theorem \ref{diffF} is proved. The second
part follows easily from Theorem \ref{1.57}(2).}  \hfill $\square$
\end{nonsec}

\begin{remark} {\rm
Since $\lim_{x,y\to 0+}D_F(x,y)=1\,,$ the constant on the right
hand side of Theorem \ref{diffF} is best possible.}
\end{remark}

\begin{question} {\rm
What is the best possible constant on the left hand side of
Theorem \ref{diffF}?}
\end{question}

\bigskip

\begin{nonsec}{\bf Proof of Theorem \ref{hRthm}.} { \rm Analogously to the proof of Theorem \ref{quotF},
we have}
$$
Q_g(x,y)=\frac{xF(x)+yF(y)}{zF(z)}\le
\frac{(x+y)F(\max\{x,y\})}{zF(\max\{x,y\})}=\frac{x+y}{z}<2.
$$
\hfill $\square$
\end{nonsec}

\bigskip

\begin{remark} \label{moregeneral} {\rm
As it is seen from the proofs, above results are valid for much
more general class of positive and monotone increasing (not
necessary differentiable) functions.

 In this sense, as the referee notes, a direct proof of the
assertion from Theorem \ref{diffF} is possible.

Namely, suppose that
$$
f_1(x):=\frac{f(x)-1}{\log(1/(1-x))}
$$
is monotone increasing for $x\in (0,1)$. Then
$$
D_f(x,y)=1+(f_1(x)-f_1(x+y-xy))\log(1/(1-x))+(f_1(y)-f_1(x+y-xy))\log(1/(1-y))
$$
$$
<1.
$$}
\end{remark}

\bigskip

\begin{nonsec}{\bf Proof of Theorem \ref{hBthm}.}
{\rm  Lemma \ref{pvlem} (1) yields
 $$
\frac{1}{B}\log(1/(1-u))\le uF(u)\le \log(1/(1-u)),
 $$
 for $u\in (0,1), cd\le 1$.

 Therefore
$$
\displaystyle\frac{x F(x)  + y F(y) }{(x+y-xy)\, F(x+y-xy) } \le
\displaystyle\frac{   \log \frac{1}{1-x} + \log
\frac{1}{1-y}}{\frac{1}{B}\log \frac{1}{(1-x)(1-y)}} = B(c,d)\,.
$$
The lower bound is proved in the same way.

Applying part (2) of Lemma \ref{pvlem}, we bound $Q_g$ similarly  in the
case $1/c+1/d\le 2$.} \hfill $\square$
\end{nonsec}

\bigskip

\begin{remark} \label{Rrmk} {\rm
From parts (1) and (3) of Lemma \ref{pvlem}, we conclude that
$$
R(c,d)>B(c,d)-1>0,
$$
for} $c,d>0, cd\le 1, (c,d)\neq (1,1)$.
\end{remark}

\begin{nonsec}{\bf Proof of Theorem \ref{Bgt1}.}
{\rm Let us write $B = B(c,d)$, $R = R(c,d)$ and
$L=\log (1/((1-x)(1-y)))>0$. By Lemma \ref{pvlem} (3) we have
  \begin{equation}\label{boundsforh}
    \frac{B-1}{B}x+\frac{1}{B}\log \frac{1}{1-x} < x F(c,d;c+d;x) < \frac{Rx}{B}+\frac{1}{B}\log \frac{1}{1-x}.
  \end{equation}
  Since $x+y < 2(x+y-xy)$ we obtain by \eqref{boundsforh}
  \begin{equation*}
    Q_g(x,y)  \le  \frac{\frac{R(x+y)}{B}+\frac{L}{B}}{\frac{B-1}{B}(x+y-xy)+\frac{L}{B}} \le \frac{2R(x+y-xy)+L}{(B-1)(x+y-xy)+L}  \le \frac{2R}{B-1}
    \,,
  \end{equation*}
  and}
  \[
    Q_g(x,y) \ge \frac{(B-1)(x+y)+L}{R(x+y-xy)+L} \ge
\frac{(B-1)(x+y-xy)+L}{R(x+y-xy)+L} \ge \frac{B-1 }{ R }\, .
  \]
\hfill $\square$
\end{nonsec}
\bigskip

\begin{nonsec}{\bf Proof of  Theorem \ref{ssthm3}.} {\rm  By the first part of Lemma \ref{pvlem}, $f$ is
monotone decreasing for $cd\le 1$.

Hence, for $0<x<y<1$ we have
$$
\frac{xF(c,d;c+d;x)}  {\log (1/ (1-x))}\ge \frac{yF(c,d;c+d;y)}
{\log (1/ (1-y))}.
$$
Putting $1-x=e^{-u}, 1-y=e^{-v}; u,v\in(0,\infty)$, we get that
the inequality
$$
\frac{(1-e^{-u})F(c,d;c+d;(1-e^{-u}))}{u}\ge
\frac{(1-e^{-v})F(c,d;c+d;(1-e^{-v}))}{v}
$$
holds whenever $0<u<v<\infty$.

This means that the function $G(t)/t$ is monotone decreasing,
where
$$G(t):=(1-e^{-t})F(c,d;c+d;(1-e^{-t}))=g(1-e^{-t})\, .$$
By Lemma \ref{mypropo}, it follows that $G$ is super-additive,
that is
$$
G(u)+G(v)\ge G(u+v),
$$
which is equivalent to
$$
g(x)+g(y)\ge g(x+y-xy),
$$
and the proof of the first part of Theorem \ref{ssthm3} is
complete.

The proof of the second part is similar. Note that the condition
$c\in (1/2,\infty), d\ge c/(2c-1)$ of Lemma \ref{pvlem} is
equivalent to the condition $1/c+1/d\le 2$ of Theorem
\ref{ssthm3}.} \hfill $\square$
\end{nonsec}

\bigskip

\begin{nonsec}{\bf Proof of Theorem \ref{ssthm4}.} {\rm (1)  The left-hand side of this inequality is a
direct consequence of part (1) of Theorem \ref{ssthm3}.

Next, from Lemma \ref{pvlem}, part (3) for $u\in (0,1), cd\le 1$, we get
$$
(B-1)u-\log(1-u)<BuF(u)< Ru-\log(1-u).
$$

Hence, by the terminology from Theorem \ref{hRthm}, we obtain
$$
BD_g(x,y)=BxF(x)+ByF(y)-BzF(z)<\log(1-z)-\log(1-x)-\log(1-y)+R(x+y)-(B-1)z
$$
$$
=(R-B+1)(x+y)+(B-1)xy < 2(R-B+1)+(B-1)=2R+1-B,
$$
since Remark \ref{Rrmk} yields $R-B+1>0$ and $B-1>0$.
\bigskip

(2) \ To prove this part we shall use Lemma \ref{pvlem}, part (4).
Because $d\ge c/(2c-1)\ge (c+1)/(3c-1)$ and $(1/2,\infty)\subset
(1/3,\infty)$, we conclude that this assertion is valid under the
condition $1/c+1/d\le 2$.

Therefore, for $u\in (0,1), 1/c+1/d\le 2, (c,d)\neq(1,1)$, we get
$$
Ru-\log(1-u)<BuF(u)< (B-1)u-\log(1-u),
$$
and, as above,
$$
BD_g(x,y)=BxF(x)+ByF(y)-BzF(z)>\log(1-z)-\log(1-x)-\log(1-y)+R(x+y)-(B-1)z
$$
$$
=(R-B+1)(x+y)+(B-1)xy > 2(R-B+1)+(B-1)=2R+1-B,
$$
since parts (1) and (4) of Lemma \ref{pvlem} give $R-B+1<0$ and $B-1<0$.

Because the right-hand inequality is a consequence of Theorem \ref{ssthm3},
part (2), the proof is complete.} \hfill $\square$
\end{nonsec}

\bigskip

\begin{nonsec}{\bf Proof of Theorem \ref{ssthm5}.} {\rm Putting
$x=y=1-\sqrt{1-t}$, we obtain $z=x+y-xy=t$. Therefore,
$$
Q_g(x,y)=2/S(t).
$$
The rest is an application of Corollary \ref{cor38}.} \hfill $\square$
\end{nonsec}

\bigskip

\subsection*{Acknowledgments}
The research of Matti Vuorinen was supported by the Academy of Finland,
Project 2600066611. The authors are indebted to the referee for his/her useful remarks.



\end{document}